%% file: FAB.tex
\documentclass[reqno]{amsart}
\pdfoutput=1 

\usepackage{general,resistance,pxfonts,ifthen,leading}
\usepackage[switch*,pagewise]{lineno}
 
\usepackage[bookmarks, colorlinks=true, pdfstartview=FitV, linkcolor=blue, citecolor=blue, urlcolor=blue]{hyperref}

\usepackage[left=1.85in,right=1.85in,top=1.7in,bottom=1.40in]{geometry}

\newcommand{\simt}{\ensuremath{\gF}\xspace}     
\newcommand{\PB}{\protect{\ensuremath{\sP}\negsp[3]}\xspace}     
\newcommand{\SG}{\protect{\ensuremath{\sS\negsp[3]\raisebox{-1pt}{\sG}\negsp[3]}}\xspace}     
\newcommand{\attr}{\ensuremath{F}\xspace}     
\newcommand{\V}{\ensuremath{\tilde V}\xspace}     
\newcommand{\G}{\ensuremath{\tilde \gG}\xspace}     
\newcommand{\Words}{\sW\xspace}     
\renewcommand{\Paths}{\gW}     
\renewcommand{\cpath}{\gw}     
\renewcommand{\attr}{\sF}     
\newcommand{\network}{\sN}     
\newcommand{\NF}{\ensuremath{\network_\attr}\xspace}     
\newcommand{\NB}{\ensuremath{\widehat{\network_\attr}}\xspace}     
\newcommand{\BF}{\ensuremath{\sB_\attr}\xspace}     
\newcommand{\vr}{\mspace{3mu}\raisebox{-3pt}{\rule{0.5pt}{12pt}}\mspace{3mu}}
\newcommand{\uKm}[1][m]{\ensuremath{u_{K^{(#1)}}}\xspace}
\newcommand{\rmax}{\ensuremath{r_{\text{max}}}\xspace}     

\numberwithin{figure}{section}
\numberwithin{equation}{section}
\numberwithin{theorem}{section}

\begin{document}
  \leading{13.5pt} 
  \linenumbers
  
\title{Self-similar fractals as boundaries of networks}

\author{Erin P. J. Pearse}
\address{University of Oklahoma, Norman, OK 73019-0315 USA}
\email{ep@ou.edu}

\thanks{\textbf{\today}. The work of EPJP was partially supported by the University of Iowa Department of Mathematics NSF VIGRE grant DMS-0602242 and by a travel grant from the Chinese University of Hong Kong.}

\begin{abstract}
  For a given pcf self-similar fractal, a certain network (weighted graph) is constructed whose ideal boundary is (homeomorphic to) the fractal. This construction is the first representation of a connected self-similar fractal as the boundary of a \emph{reversible} Markov chain (i.e., a simple random walk on a network). The boundary construction is effected using certain functions of finite energy which behave like bump functions on the boundary. The random walk is shown to converge to the boundary almost surely, with respect to the standard measure on its trajectory space.
\end{abstract}




\maketitle

\setcounter{tocdepth}{1}
{\small \tableofcontents}

\allowdisplaybreaks

\input{introduction}

\input{basic-terms}

\input{network-construction}

\input{random-walk-measures}

\input{separating-points}

\input{convergence}

\input{identification}

\subsection*{Acknowledgements}

The author is grateful for stimulating comments, helpful advice, and valuable references from Jun Kigami, Ka-Sing Lau, Sze-Man Ngai, Wolfgang Woess, and the students and colleagues who have endured talks on this material and raised fruitful questions. In particular, I am grateful to Ka-Sing Lau for the invitation to visit him at the Chinese University of Hong Kong and to both Ka-Sing Lau and Sze-Man Ngai for explaining several things to me while I was there; especially in regard to Martin boundaries. I would also like to thank Jun Kigami for answering several questions and for bringing the remarkable paper \cite{AnconaLyonsPeres} to my attention. I was greatly motivated by conversations with Jun, and by \cite{Kig01,Kig09b}. Finally, I would like to thank Palle Jorgensen for extensive discussions and advice, and for being so generous with his time during my postdoc at U. Iowa.

\bibliographystyle{alpha}
\bibliography{networks}

\end{document}

%% file: introduction.tex

\section{Introduction}

In this paper, we construct a network that has a given self-similar fractal as its boundary. The intended application of such a network is as a means for further studying the theory of analysis on fractals as developed by Kigami \cite{Kig01} and others (in particular, by Strichartz, Teplyaev, and their associates, see \cite{Str06} and its references). The motivation for understanding a fractal as a boundary is as follows: one can obtain a Laplace operator on the fractal as the trace of a discrete Laplacian on the network, and it is much easier to work with the latter. More precisely, one can obtain a Dirichlet form on the boundary as the trace of a Dirichlet form on the network, and the associated Laplacian can then be understood via the usual correspondence given by Kato's theorem; cf.~\cite[Thm.~1.3.1]{FOT}. Whether or not this Dirichlet form corresponds to a diffusion is another question; it is more likely that the present construction corresponds to a jump process; see.
Nonetheless, this opens up a new avenue to the study of analysis on fractals; one which may generalize to non-post-critically finite examples (see \cite[Def.~1.3.12]{Kig01} for the definition of p.c.f.). Using different methods, Kigami has made significant inroads in this direction for self-similar fractals which are totally disconnected. See \cite{Kig09b} for explicit formulas of the induced energy form on the boundary, estimates for the kernel of the associated jump process, and a host of other interesting results.

Consider a self-similar fractal set \attr defined by an iterated function system $\{\simt_0, \simt_1, \dots, \simt_J\}$.%
  \footnote{The fractal is assumed to be post-critically finite (pcf) and satisfy a regularity condition; see Axiom~\ref{axm:pcf} in \S\ref{sec:the-network}.} We construct a connected graph \NF which has the fractal \attr as its boundary, in the sense of \cite[\S7]{WoessII}: 
\begin{defn}\label{def:compactification}
  a \emph{compactification} of a countably infinite set $X$ is a compact Hausdorff space $\widehat X$ containing $X$ and satisfying
  \begin{enumerate}[(i)]
    \item the set $X$ is dense in $\widehat X$, and
    \item in the induced (subspace) topology, $X$ is a discrete subset of $\widehat X$, 
  \end{enumerate}
  and the \emph{boundary} of $X$ is defined to be $\sB := \widehat X \less X$.  
\end{defn}
 
Such a compactification can be constructed in terms of a family of functions so as to have certain desirable characteristics. The following result is \cite[Thm.~7.13]{WoessII}.
\begin{theorem}\label{thm:compactification}
  If \xF is a countable family of bounded functions on $X$, then there exists a unique (up to the appropriate notion of equivalence) compactification $\widehat X = \widehat X_\xF$ of $X$ such that 
  \begin{enumerate}[(a)]
    \item every function $f \in \xF$ extends to a continuous function on $\widehat X$, and
    \item \xF separates boundary points: for distinct $\gx, \gh \in \sB$, there is $f \in \xF$ with $f(\gx) \neq f(\gh)$.
  \end{enumerate}
\end{theorem}

\begin{theorem}\label{thm:compactification2}
  The compactification $\widehat X_\xF$ of Theorem~\ref{thm:compactification} is equivalently characterized as $\widehat X_\xF = X \cup (\gW_\iy / \sim)$ where
  \linenopax
  \begin{align}
    \gW_\iy &= \{\gw=(x_n) \in X^\bN \suth x_n \to \iy \text{ and } f(x_n) \text{ converges for all } f \in \xF\}, 
    \label{eqn:gW-inf} \\
    \qq\text{and } 
    (x_n) &\sim (y_n) \iff \lim f(x_n) = \lim f(y_n), \text{ for every } f \in \xF.
    \label{eqn:sim}
  \end{align}
\end{theorem}
The network \NF is compactified in Definition~\ref{def:boundary} according to Theorem~\ref{thm:compactification2}, and Theorem~\ref{thm:Convergence-to-the-boundary} shows that the random walk on \NF converges almost surely (with respect to the natural path-space measure) to a point of $\BF = \NB \less \NF$, in the topology of \NB. In Theorem~\ref{thm:B=F}, we use Theorem~\ref{thm:compactification} to verify that the boundary so obtained coincides with the fractal \attr. See Example~\ref{exm:SN} (Figure~\ref{fig:sierpinski-bd-fn}) and Example~\ref{exm:C} (Figure~\ref{fig:cantor-network}) for the case of the Sierpinski gasket and the Cantor set.
  
  For construction outlined above, \xF will be a subfamily of the functions of finite energy on \NF; these are constructed in \S\ref{sec:separating-points}. For the classical Martin boundary, one uses $\xF = \{K(x,\cdot)\}_{x \in X}$, where the Martin kernel $K(x,y)$ is defined in terms of the Green kernel (see \eqref{eqn:Green-function}) via
  \linenopax
  \begin{align}\label{eqn:Martin-kernel}
    K(x,y) = \frac{G(x,y)}{G(o,y)}, \qq x,y \in X,
  \end{align}
for some fixed reference vertex $o \in X$. This definition makes each $K(x,\cdot)$ a bounded and superharmonic function, so that $K(x,Z_n)$ is a supermartingale with respect to the random walk $(Z_n)_{n=0}^\iy$ on $X$. From this supermartingale, one can compute a crossings estimate, and hence obtain results about the convergence of $Z_n$ and $K(x,Z_n)$ necessary for making the notion of boundary rigourous, and solving the Dirichlet problem at infinity. 

Martin kernels are well-suited to studying positive harmonic functions, but not as well-adapted for studying functions of finite energy. We must replace the Martin kernels with some family of finite-energy functions which retain the convergence properties of supermartingales described above. A key observation in \cite{AnconaLyonsPeres} is that if $f$ is a function of finite Dirichlet energy ($\energy(f)<\iy$ in Definition~\ref{def:energy}), then one can find a superharmonic function $h$ with $\energy(h)<\energy(f)$ and $h(x) \geq |f(x)|$, for all $x \in X$. This $h$ gives rise to a supermartingale from which one can obtain a crossings estimate, and we apply this idea to obtain a boundary theory for \NF. 

In \cite[\S7]{WoessII}, Woess remarks that the measurable structure of the Poisson boundary of a network makes it the ``right'' model for distinguishable limit points at infinity which the random walk can attain. The boundary measure space $\PB = (\attr, \gn_o)$ is the Poisson boundary of \NF, where $\gn_o$ is a certain probability measure related to the random walk on \NF.
This would yield 
a unique integral representation of the harmonic functions of finite energy:
  \linenopax
  \begin{align}\label{eqn:integral-representation-intro}
    h_\gf(x) = \int_\PB \gf \, d\gn_x = \Ex_x(\gf(X_\iy)),
    \qq\text{for all } x \in \NF,
  \end{align}
where $\gf \in L^\iy(\PB)$ is prescribed, and $\gn_x$ is a certain Radon-Nikodym derivative of $\gn_o$.
This allows one to obtain an energy form on the fractal as the trace of the energy form on the network. More precisely, for $\gf \in L^\iy(\PB)$, formula \eqref{eqn:integral-representation-intro} gives the harmonic extension $h = h_\gf$ of \gf to \NF, whence one defines the trace energy form $\energy_\PB$ via
  \linenopax
  \begin{align}\label{eqn:trace-form}
    \energy_\PB(\gf) := \energy(h_\gf),
    \qq \dom \energy_\PB = \{\gf \in L^\iy(\PB) \suth h_\gf \in \dom \energy\}.
  \end{align}
In this way, one obtains a (Dirichlet) energy form $\energy_\PB$ on the fractal. This approach has already been successful for studying trees with totally disconnected fractal boundary: see \cite{Kig09b}, where an extended discussion of \eqref{eqn:integral-representation-intro}--\eqref{eqn:trace-form} is given in the context of the Douglas integral.

For the network $\NF = (V, E \cup F)$ constructed below, the vertex set is denoted $V$, and $E$ and $F$ are sets of ``horizontal edges'' and ``vertical edges'', respectively, following Kaimanovich's terms in \cite{Kaimanovich}. However, this is not the same as the graphs described by Kaimanovich. In particular, the Sierpinski \emph{network} discussed here is different from the Sierpinski \emph{graph} in \cite{Kaimanovich}, where the author shows how the Sierpinski gasket can be interpreted as the hyperbolic boundary of a graph, in the sense of Gromov. Also, the present construction is rather different from the nonreversible Markov processes (corresponding to \emph{directed} graphs) considered by Denker \& Sato \cite{DenkerSatoI, DenkerSatoII} and also by Ju, Lau, and Wang \cite{JuLauWang}. (See also \cite{LauWang} for a representation of a pcf self-similar fractal as the hyperbolic boundary of a graph.) The paper \cite{Cartier} provides a very readable account of harmonic functions and boundaries specific to trees.

%% file: basic-terms.tex

\section{Basic terms}
\label{sec:basic-terms}

\begin{defn}\label{def:network}
  A \emph{network} is a connected graph $(G,\cond)$, where $G = (V,E)$ is a connected simple graph with vertex set $V$ and edge set $E$, and \cond is a \emph{conductance function} satisfying
  \begin{enumerate}[(i)]
    \item $c_{xy}>0$ iff $x$ and $y$ are adjacent vertices (denoted $x \nbr y$), and
    \item $\cond_{xy} = \cond_{yx} \in [0,\iy)$.
  \end{enumerate}
  for all vertices $x$ and $y$ (which is henceforth denoted $x,y \in G$). We write $\cond(x) := \sum_{y \nbr x} \cond_{xy}$ and require $\cond(x) < \iy$, for each fixed $x$. 
  In this definition, \emph{connected} means that for any $x,y \in G$, there is a \emph{path} connecting $x$ to $y$, i.e., there is a finite sequence $\{x_i\}_{i=0}^n$ with $x=x_0$, $y=x_n$, and $\cond_{x_{i-1} x_i} > 0$, $i=1,\dots,n$. 
\end{defn}

\begin{defn}\label{def:Laplacian}
  The \emph{Laplacian} on G is the linear difference operator which acts on a function $v:G \to \bR$ by
  \linenopax
  \begin{equation}\label{eqn:Laplacian}
    (\Lap v)(x) := \sum_{y \nbr x} \cond_{xy}(v(x)-v(y)).
  \end{equation}
  A function $v:G \to \bR$ is called \emph{harmonic} iff $\Lap v \equiv 0$. The domain of \Lap is specified in Definition~\ref{def:vx}.
\end{defn}


\begin{defn}\label{def:energy}
  The \emph{energy} of a function $u:G \to \bR$ is given by the Dirichlet form
  \linenopax
  \begin{align}\label{eqn:energy-form}
    \energy(u)
    := \; \frac12 \sum_{x,y \in G} \cond_{xy}(u(x)-u(y))^2.
  \end{align}
  Each summand in \eqref{eqn:energy-form} is nonzero iff $x \nbr y$, so this is really a sum over the edges in the network, whence the initial factor of $\frac12$. 
  The \emph{domain of the energy} is
  \linenopax
  \begin{equation}\label{eqn:energy-domain}
    \dom \energy = \{u:G \to \bR \suth \energy(u)<\iy\}.
  \end{equation}
\end{defn}

\begin{defn}\label{def:HE}
  The space $\HE := \dom \energy / \bR \one$ is a Hilbert space under the inner product
  \linenopax
  \begin{align}\label{eqn:energy-inner-product}
    \la u, v \ra_\energy 
    := \frac12 \sum_{x,y \in G} \cond_{xy}(u(x)-u(y))(v(x)-v(y)),
  \end{align}
  and with the corresponding norm $\|u\|_\energy = \sqrt{\la u, u \ra_\energy} =\energy(u)^{1/2}$; see \cite{DGG,ERM,Kig01,Kig03}. The space \HE will be called the \emph{energy Hilbert space} and \eqref{eqn:energy-inner-product} the \emph{energy product}. 
\end{defn}
  
Note that since \eqref{eqn:energy-form} is a sum of nonnegative terms, its convergence is independent of rearrangements. Consequently, the sum in \eqref{eqn:energy-inner-product} is well-defined by the Cauchy-Schwarz inequality.
 
\begin{defn}\label{def:vx}\label{def:energy-kernel}
  Let $v_x$ be defined to be the unique element of \HE for which
  \linenopax
  \begin{equation}\label{eqn:v_x}
    \la v_x, u\ra_\energy = u(x)-u(o),
    \qq \text{for every } u \in \HE.
  \end{equation}
  The collection $\{v_x\}_{x \in G}$ forms a reproducing kernel for \HE (\cite[Cor.~2.7]{DGG}); it is called the \emph{energy kernel} and \eqref{eqn:v_x} shows its span is dense in \HE. Define $\dom \Lap = \spn\{v_x\}_{x \in G}$.
\end{defn}

\begin{defn}\label{def:R(x)}
  Denote the \emph{(free) effective resistance} between $x$ and $y$ by 
  \linenopax
  \begin{align}\label{eqn:R(x)}
    R(x,y) = R^F(x,y) := \energy(v_x - v_y).
  \end{align}
  This quantity represents the voltage drop measured when one unit of current is passed into the network at $x$ and removed at $y$, and the equalities in \eqref{eqn:R(x)} are proved in \cite{ERM} and elsewhere in the literature; see \cite{Lyons,Kig03} for different formulations. 
\end{defn}

The following results appear in \cite[Lem.~2.23]{DGG} and \cite[Lem.~6.9]{bdG} and will be useful in the sequel; for further details, see \cite{DGG, ERM, Multipliers, SRAMO, bdG, RBIN, RANR, Friedrichs, Interpolation, LPS, OTERN}.

\begin{lemma}
  \label{thm:energy-kernel-is-real}\label{thm:monopoles-and-dipoles-are-bounded}
  Every $v_x$ is \bR-valued, with $v_x(y) - v_x(o) >0$ for all $y \neq o$. 
  Every $v_x$ is bounded with $\|v_x\|_\iy := \sup_{y \in G}|v_x(y)-v_x(o)| = v_x(x)-v_x(o) = R(x,o) < \iy$.
\end{lemma}

\begin{remark}\label{rem:v_x-reps}
  It will be useful (especially in \S\ref{sec:convergence}) to discuss energy kernel elements $v_x$ in the context of functions on the network. In this case, we abuse notation and also write $v_x$ for the function on \NF which is the representative of $v_x$ satisfying $v_x(o)=0$. In such a context, $v_x(y)=v_y(x)$ because \eqref{eqn:v_x} implies
  \linenopax
  \begin{align*}
    v_y(x) - v_y(o) 
    = \la v_x, v_y\ra_\energy 
    = \la v_x, v_y\ra_\energy
    = v_x(y) - v_x(o),
  \end{align*}
  where the central equality follows from Lemma~\ref{thm:energy-kernel-is-real}. It is shown in \cite{DGG} that 
  \linenopax
  \begin{equation}\label{eqn:dipole}
    \Lap v_x = \gd_x-\gd_o,
  \end{equation}
  where $\gd_x$ denotes the characteristic function of $x$.
\end{remark}

%% file: network-construction.tex

\section{The self-similar fractal and the network}
\label{sec:the-network}

The class of self-similar fractal sets \attr under investigation is defined in \S\ref{sec:fractals} and the corresponding network \NF is constructed in \S\ref{sec:network-construction}. Its vertices form a set $V \ci \attr \times \bZ_+$, where $\bZ_+ = \{0,1,2,\dots\}$, and its edges are defined according to the iterated function system defining \attr. 

\subsection{Self-similar fractals}
\label{sec:fractals}

Let $\{\simt_0, \simt_1, \dots, \simt_J\}$ be an iterated function system (IFS) of finitely many contractive similarity mappings on \bRd where each $\simt_j$ is a composition of a rigid motion of \bRd and a homothety with scaling factor $r_j \in (0,1)$, for $j=0,1,\dots,J$. Also, each $\simt_j$ has unique fixed point $q_j$. Let $\V_0 = \{q_0,q_1,\dots,q_J\}$ and define the action of \simt on the space of compact subsets of \bRd by 
\linenopax
\begin{align}\label{eqn:simt}
  \simt\left(A\right) := \bigcup_{j=0}^J \simt_j \left(A\right).
\end{align}

  The (fractal) \emph{attractor} \attr is the unique nonempty compact subset of \bRd satisfying the fixed-point equation $\attr = \simt(\attr)$; cf.~\cite{Hut81}. 

\begin{defn}\label{def:m-cell}
  If $\Words_m := \{0,1,\dots,J\}^m$ is the set of \emph{words} $w = w_1 w_2 \dots w_m$ of length $m$, then we write the length of $w \in \Words_m$ as $|w| = m$, and use the following shorthand notation for a composition of the mappings $\simt_j$:
  \linenopax
  \begin{align}\label{eqn:word-composition}
    \simt_w(x) := \simt_{w_1} \comp \simt_{w_2} \comp \dots \comp \simt_{w_m} (x).
  \end{align}
  Let $\Words_0 = \{\es\}$, where $\es$ is the unique word of length $0$, and let $\simt_\es = \id$ be the identity mapping. For a word $w \in \Words_m$, the set $\simt_w(\attr)$ is called a \emph{$m$-cell}. We denote the set of all finite words by $\Words_\star := \bigcup_{m=0}^\iy \Words_m$ and the set of all infinite words by \Words.
\end{defn}

There is a continuous surjection 
  \linenopax
  \begin{align}\label{eqn:quotient}
    \gp:\Words \to \attr
    \qq\text{given by}\qq
    \gp(w) = \bigcap_{n=1}^\iy \simt_{w|n}(\attr)
  \end{align}
where $w|n := w_1 w_2 \dots w_n$ is the truncation of $w$. Thus \attr is a quotient of \Words by the equivalence relation $w \sim_\attr w'$ iff $\gp(w) = \gp(w')$. Points $\gx=\gp(w)$ of the fractal \attr can be approximated most easily%
\footnote{In light of \eqref{eqn:quotient}, it is clear that $\lim_{n \to \iy}\simt_{w|n}(q) = \gx$ for any $q \in \attr$, but it is convenient (and sufficient) to consider $q \in \V_0$.}
  by sequences of rational points (see Definition~\ref{def:junction-point}): 
  \linenopax
  \begin{align}\label{eqn:rational-approximants}
    (\gx_n)_{n=1}^\iy, 
    \qq\text{where}\qq 
    \gx_n := \simt_{w|n}(q),
    \text{ for some fixed } q \in \V_0.  
  \end{align}
This leads to the discretization described in Definition~\ref{def:Vk-and-Gk} and is also the idea behind the network construction in \S\ref{sec:network-construction}.

\begin{figure}
  \scalebox{0.80}{\includegraphics{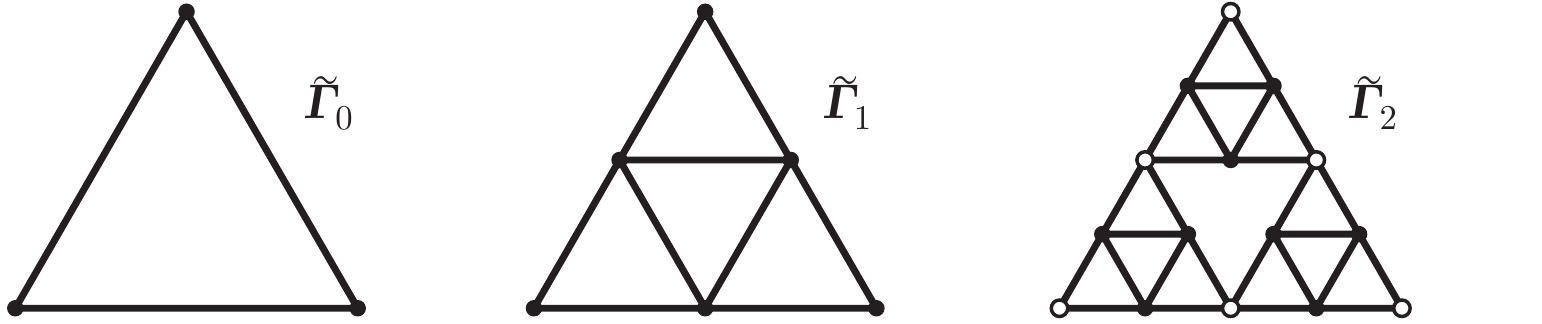}}
  \centering
  \caption{\captionsize The discrete approximant $\G_k$ of the Sierpinski gasket is isomorphic to the corresponding subnetwork $\gG_k$ (at level $k$) of the Sierpinski network $\network_{\negsp[4]\sS\sG}$ constructed in \S\ref{sec:network-construction}; see Figure~\ref{fig:sierpinski-bd-fn}. The white vertices in $\G_2$ form the set $V_1^{(2)}$ as in \eqref{eqn:G_k^m}.} 
  \label{fig:sierpinski-network-h-edges}
\end{figure}

\begin{defn}[Discrete prefractal approximants]
  \label{def:Vk-and-Gk}
  We define 
  \linenopax
  \begin{align}\label{eqn:Vtilde_k}
    \V_k := \simt(\V_{k-1}) = \simt^k(\V_0),
  \end{align}
  and refer to \eqref{eqn:Vtilde_k} as the set of points arising in the \kth stage of the construction of \attr, where \simt is given by \eqref{eqn:simt}, and $\simt^k$ refers to the $k$-fold application of \simt. 
  Let $\G_0$ be the complete network on the fixed points $\V_0$, and fix some choice of edge weights $c_{xy}=c_{yx}$, for $x,y \in \V_0$. Finally, let $\G_k = \simt(\G_{k-1}) = \simt^k(\G_0)$ be the graph with vertex set $\V_k$ and with edges defined by 
  \linenopax
  \begin{align}\label{eqn:cond}
    \cond_{\simt_j(x) \simt_j(y)} = \frac{\cond_{xy}}{r_j}, \qq \text{for each } j=0,1,\dots,J,
  \end{align}
  and for some choice of renormalization factors $r_j > 0$, for $j=0,1,\dots,J$. 
  The sets $\G_k$ are illustrated for the Sierpinski gasket in Figure~\ref{fig:sierpinski-network-h-edges}.
\end{defn}

We use the tilde notation $\V$ (and $\V_k$, etc) to refer to vertices of the fractal \attr and the non-tilde version $V$ (and $V_k$, etc) to refer to vertices of the network \NF, in \S\ref{sec:network-construction}. 
We will require the existence of a regular self-similar energy form $\energy_\attr$ on \attr, as described in \cite[\S3.1]{Kig01} or \cite[\S1.4 and \S4.2]{Str06} in detail, and briefly reviewed here. 

\begin{defn}[Self-similar energy form]\label{def:Self-similar-energy-form}
  For any function $u:\G_m \to \bR$, the energy $\energy_m(u)$ is defined via \eqref{eqn:energy-form} on the finite graph $\G_m$. Since $\V_{m-1} \ci \V_m$, one can consider the restriction $u|_{\G_{m-1}}$ by considering the values of $u$ on $\V_{m-1}$, and computing $\energy_{m-1}(u|_{\G_{m-1}})$ via \eqref{eqn:energy-form} on the finite graph $\G_{m-1}$. Note that the edge sets of $\G_m$ and $\G_{m-1}$ are different, and in particular, the respective conductances are defined by \eqref{eqn:cond}.
  The fractal \attr admits a self-similar energy form iff the renormalization factors $r_j$ and the conductances $c_{xy}$ of Definition~\ref{def:Vk-and-Gk} can be chosen so that
  \linenopax
  \begin{align}\label{eqn:m-energy}
    \energy_m(u) = \energy_{m-1}(u|_{\G_{m-1}})
    \qq\text{and}\qq
    \energy_m(u) = \sum_{j=0}^J r_j^{-1} \energy_{m-1}(u \comp \simt_j).
  \end{align}
  When \eqref{eqn:m-energy} is satisfied, one can define
  \linenopax
  \begin{align}\label{eqn:energy-def}
    \energy_\attr(u) = \lim_{m \to \iy} \energy_m(u),
  \end{align}
  and be assured that the limit exists in $[0,\iy]$ and has the self-similar property
  \linenopax
  \begin{align}\label{eqn:self-similar-energy}
    \energy_\attr(u) = \sum_{j=0}^J r_j^{-1} \energy_ \attr(u \comp \simt_j).
  \end{align}
  The domain of the quadratic form \eqref{eqn:energy-def} is $\dom \energy_\attr := \{u:\attr\to \bR \suth \energy_\attr(u) < \iy\}$.
\end{defn}

\begin{axm}[Regularity]\label{axm:pcf}
  The self-similar set \attr is required throughout to be \emph{post-critically finite} (pcf), as in \cite[Def.~1.3.13]{Kig01}, and to admit a \emph{regular} harmonic structure, as in \cite[Def.~3.1.2]{Kig01}. This regularity condition means that one can find a self-similar energy form $\energy_\attr$ for which
  \linenopax
  \begin{align}\label{eqn:r_j<1}
    r_j < 1, \qq\text{for each}\qq j=0,1,\dots,J. 
  \end{align}
\end{axm}

\begin{remark}\label{rem:axiom}
  Axiom~\ref{axm:pcf} is a technical tool, and it is possible that it may be removed in the future. Regularity is only used in the proof of Lemma~\ref{thm:harmext-is-finite-energy} to show that certain functions have finite energy. Axiom~\ref{axm:pcf} is satisfied in most cases, no general conditions are known that guarantee it; see the discussion following \cite[Prop.~3.1.3]{Kig01} or \cite[p.98--99]{Str06}.
\end{remark}

\begin{exm}[The Sierpinski gasket]\label{exm:SG}
  The Sierpinski gasket \SG is defined by an iterated function system of three mappings $\simt_j = \frac{x-q_j}2 + q_j$, $j=0,1,2$, where $\{q_0,q_1,q_2\}$ are any three noncollinear points in \bRt. Note that $q_j$ is the unique fixed point of $\simt_j$. The standard choice of conductances for \SG is to take $c_{xy} \equiv 1$ for $x,y in \V_0$ (and hence for all $x,y \in \V$, by \eqref{eqn:cond}), and the standard choice of renormalization factors is to take $r_0=r_1=r_2=\frac35$. 
\end{exm}

\subsection{Construction of the network}
\label{sec:network-construction}

 Let us use a superscript for subsets of $\attr \times \bZ_+$ to indicate the ``vertical coordinate'', as in the following definitions.

\begin{figure}
  \scalebox{0.80}{\includegraphics{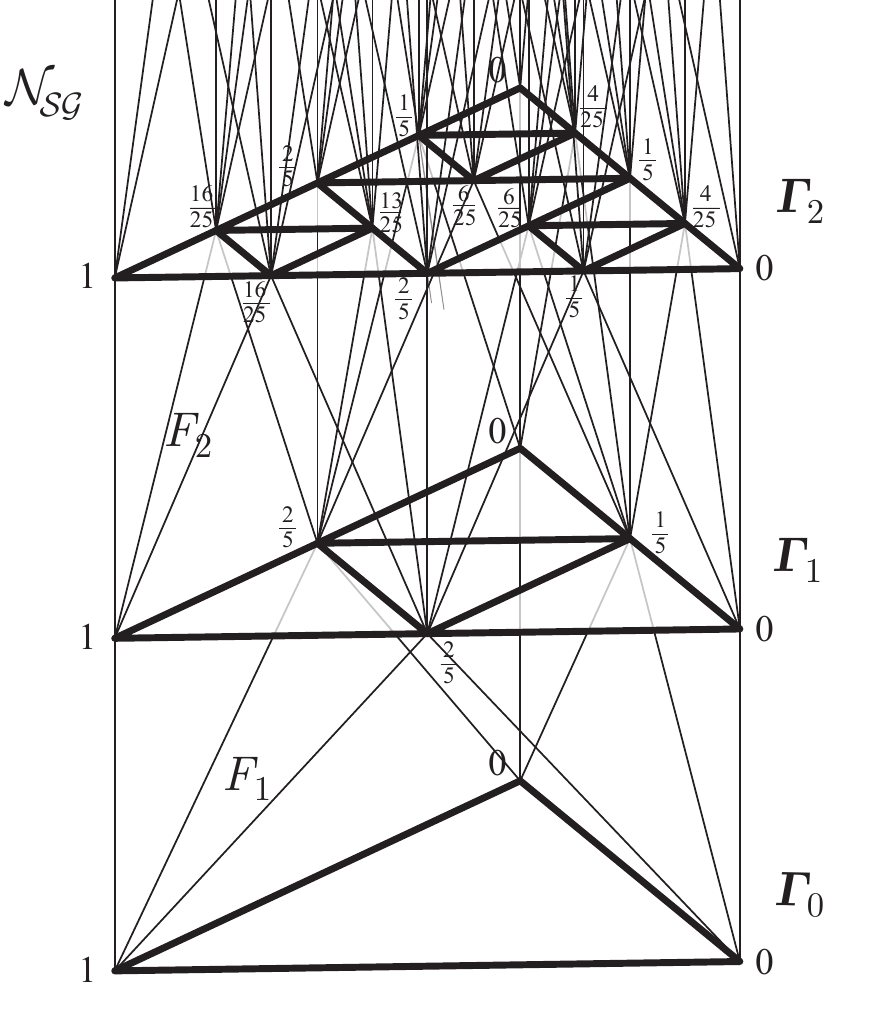}}
  \centering
  \caption{\captionsize A depiction of the Sierpinski network $\network_{\negsp[4]\sS\sG}$ corresponding to $\attr=\sS\sG$ as in Example~\ref{exm:SG}; see Definition~\ref{def:Sierpinski-network}. Note the $k$-level subgraph $\gG_k$ of $\network_{\negsp[4]\sS\sG}$ is isomorphic to $\tilde \gG_k = \simt^k(\G_0)$ as in Figure~\ref{fig:sierpinski-network-h-edges}, and also how the vertical edge sets $F_k$ connects $\gG_{k-1}$ to $\gG_k$; see Figure~\ref{fig:sierpinski-network-v-edges}. Numbers labeling the vertices are values of the function $u$ constructed in Lemma~\ref{thm:u-nonvanishing-at-infinity}.}
  \label{fig:sierpinski-bd-fn}
\end{figure}

\begin{defn}\label{def:V}
  For $\V_k \ci \attr$: 
\linenopax
  \begin{align}\label{eqn:G_k^m}
  V_k^{(m)} = \{x=(\gx,m) \in \attr \times \bZ_+ \suth \gx \in \V_k\} 
  = \V_k \times \{m\},
\end{align}
and let 
\linenopax
  \begin{align*}
  V_k := V_k^{(k)} = \V_k \times \{k\} 
  \qq\text{and}\qq
  \gG_k := \G_k \times \{k\}
\end{align*}
be the vertices of \NF at level $k$, and the subnetwork of \NF at level $k$, respectively. 
The vertex set of the network \NF is $V := \bigcup_{k=0}^\iy V_k$; see Figure~\ref{fig:sierpinski-bd-fn}. 
\end{defn}

\begin{remark}\label{rem:subgraphs}
   Note that $\gG_{k-1}$ is not a subgraph of $\gG_k$; rather, these are disjoint subgraphs of \NF. 
\end{remark}

\begin{defn}\label{def:junction-point}
  A \emph{rational point} of a pcf self-similar set is a point \gx for which one can find $q \in \V_0$ and $w \in \Words_m$ for some (finite) $m$, such that $\gx = \simt_w(q)$. The \emph{generation} of a rational point \gx is the first time it appears, i.e., the smallest $m$ for which such a representation can be found:
  \linenopax
  \begin{align}\label{eqn:generation}
    m_\gx := \min \{m \geq 0 \suth \simt_w(q)=\gx \text{ for } w \in \Words_m, q \in V_0\}.
  \end{align}
  Thus, $\V_k$ is the set of rational points of generation $k$, and to each vertex in $V_k$ there corresponds a rational point.
  
  A \emph{junction point} is rational point which has more than one such representation; for any pcf self-similar set there is an $N \in \bN$ such that any junction point \gx has at most $N$ such representations. 
\end{defn}

\begin{defn}[Horizontal edges]
  \label{def:horizontal-edge}
 Let $E_k$ denote the set of edges of the graph $\gG_k$. The set of \emph{horizontal edges} of \NF is $E = \bigcup_{k=1}^\iy E_k$, and we say that elements of $E_k$ are the horizontal edges at level $k$; see Figure~\ref{fig:sierpinski-network-h-edges}. The conductances of the horizontal edges are determined as in Definition~\ref{def:Vk-and-Gk}.
\end{defn}

\begin{figure}
  \scalebox{1.0}{\includegraphics{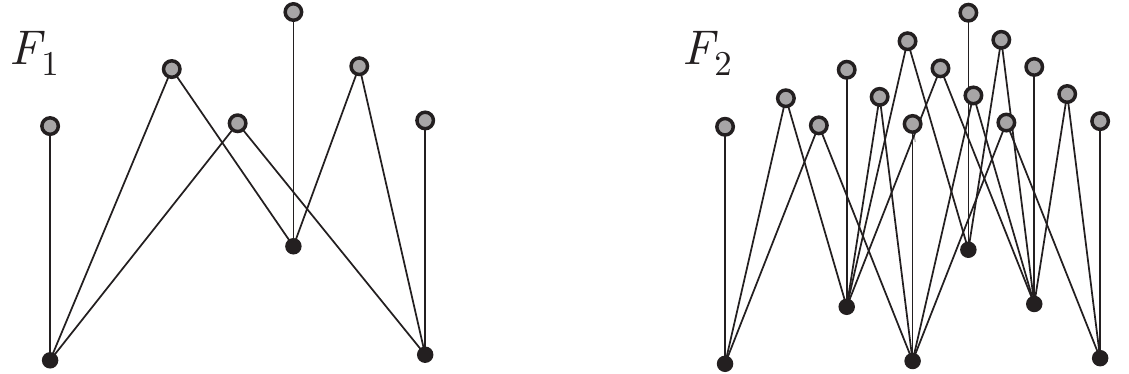}}
  \centering
  \caption{\captionsize The vertical edges of the Sierpinski network $\network_{\negsp[4]\sS\sG}$. For each $F_k$, the vertices with grey centers lie in $\gG_k$ and the vertices with black centers lie in $\gG_{k-1}$.}
  \label{fig:sierpinski-network-v-edges}
\end{figure}\vspace{-1.5ex}

\begin{defn}[Vertical edge]
  \label{def:vertical-edge}
  For each $i,j=0,1,\dots,J$ and $w \in \Words$, the network \NF contains a \emph{vertical edge} connecting $(\simt_w(q_i),k-1) \in V_{k-1}$ to  $(\simt_{wj} (q_i),k) \in V_{k}$. 
  \linenopax
  \begin{align}\label{eqn:Fk}
    F_k := \{ [(\simt_w(q_i),k-1), (\simt_{wj} (q_i),k)] \suth i,j \in \{0,1,\dots,J\}, w \in \Words_{k-1}\},
  \end{align}
  is called the set of \emph{vertical edges at level $k$}, and the set of all vertical edges is $F = \bigcup_{k=1}^\iy F_k$. See Figure~\ref{fig:sierpinski-bd-fn} and Figure~\ref{fig:sierpinski-network-v-edges}.
  
  For the conductances of the vertical edges, choose a system of probability weights to associate to the mappings in the IFS defining \attr: 
  \linenopax
  \begin{align}\label{eqn:probability-weights}
    \{\gm_0,\gm_1,\dots,\gm_J\}, 
    \q\text{where}\q 
    \gm_j > 0 \text{ for } j=0,1,\dots,J
    \q\text{and}\q 
    \sum_{j=0}^J \gm_j = 1.
  \end{align}
  Now define the vertical conductances by 
  \linenopax
  \begin{align}\label{eqn:vertical-edge-weights}
    c_{(\simt_w(q_i),k-1),(\simt_{wj} (q_i),k)} := \gm_j, 
  \end{align}
  for each vertical edge $[(\simt_w(q_i),k-1), (\simt_{wj} (q_i),k)]$, as in \eqref{eqn:Fk}.
\end{defn}

\begin{remark}\label{rem:vertical-normalization}
  In this paper, we will always choose uniform probability weights. Consequently, we may renormalize and take $c_{xy}=1$ for any edge $[x,y] \in F$, to avoid dealing with the extraneous factor of $\frac1{J+1}$.
  
  It is shown in \cite{Hut81} that for an IFS of the type described in \S\ref{sec:fractals}, there is actually an invariant \emph{measure} with support \attr, uniquely determined by a choice of probabilities $\{\gm_0,\gm_1,\dots,\gm_J\}$, where $\gm_j > 0$ for $j=0,1,\dots,J$ and $\sum_{j=0}^J \gm_j = 1$. Here, each $\gm_j$ is associated to $\simt_j$, and one has a self-similar measure satisfying
  \linenopax
  \begin{align}\label{eqn:self-similar-meas}
    \gm(A) = \sum_{j=0}^J \gm_j \cdot \gm(\simt_j^{-1} A).
  \end{align}
  See \cite{Hut81} and \cite[\S1.2]{Str06}. Also, \cite{Kig01,Kig03} describe the  harmonics structures on self-similar sets corresponding to different weightings. These generalization will be considered in a forthcoming paper.
\end{remark}

\begin{defn}\label{def:Sierpinski-network}
  The network \NF is the network with vertices $V$ as in Definition~\ref{def:V} and edges $E \cup F$ as described in Definition~\ref{def:horizontal-edge} and Definition~\ref{def:vertical-edge}. 
\end{defn}

\begin{remark}[Nonunit conductances]
  \label{rem:c=1}
  The reader may wonder why the networks \NF have all vertical edge weights equal to $1$, even though the framework in \S\ref{sec:basic-terms} allows for conductance functions $c_{xy}$ which may vary. The primary reason for working in this generality is that \cite{Hut81} shows that for an iterated function system of the type described in \S\ref{sec:fractals}, there is actually an invariant \emph{measure} with support \attr, uniquely determined by a choice of probabilities $(\gm_0,\gm_1,\dots,\gm_J)$, where $\gm_j > 0$ for $j=0,1,\dots,J$ and $\sum_{j=0}^J \gm_j = 1$. Here, each $\gm_j$ is associated to $\simt_j$, and one has a self-similar measure satisfying
  \linenopax
  \begin{align}\label{eqn:self-similar-meas}
    \gm(A) = \sum_{j=0}^J \gm_j \cdot \gm(\simt_j^{-1} A).
  \end{align}
  See \cite{Hut81} and \cite[\S1.2]{Str06}. Also, \cite{Kig01,Kig03} describe the  harmonics structures on self-similar sets corresponding to different weightings. These generalization will be considered in a forthcoming paper.
\end{remark}

\begin{exm}[Sierpinski network]\label{exm:SN}
  For the case $\attr = \SG$ of Example~\ref{exm:SG}, the Sierpinski network $\network_{\negsp[4]\sS\sG}$ is depicted in Figure~\ref{fig:sierpinski-bd-fn}. Here, $\gG_k$ contains $\frac{1}{2}(3+3^k)$ vertices and $3^k$ edges, and there are $\frac{1}{2}(3+5\cdot3^k)$ edges in $F_k$ (from $\gG_{k-1}$ to $\gG_{k}$).
\end{exm}

\begin{figure}
  \scalebox{1.0}{\includegraphics{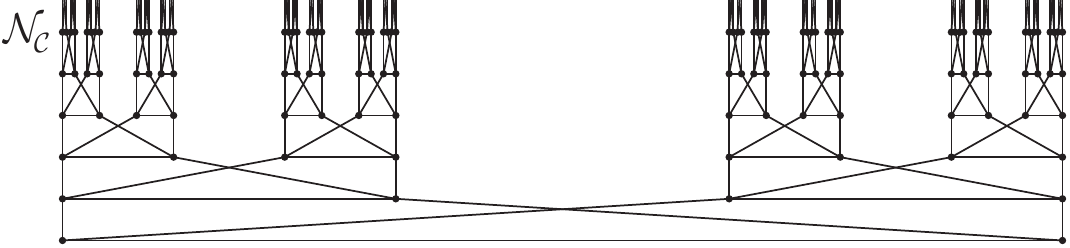}}
  \centering
  \caption{\captionsize The Cantor network $\network_{\sC}$ of Example~\ref{exm:C}. The discrete approximant $\G_k$ of the Cantor network consists of $2^k$ disjoint subgraphs, each containing two vertices and a single edge.} 
  \label{fig:cantor-network}
\end{figure}\vspace{-1.5ex}

\begin{exm}[Cantor network]\label{exm:C}
  The (standard, ternary) Cantor set \sC is defined by an iterated function system of two mappings $\simt_0(x) = \frac{x}2$ and $\simt_1(x)= \frac{x+1}2$. For the case $\attr = \sC$, the Cantor network $\network_{\sC}$ is depicted in Figure~\ref{fig:cantor-network}.
\end{exm}

\begin{defn}[Directly above and directly below]
  \label{def:directly-below}
  We say that $y = (y_1,y_2) \in V_l$ is \emph{directly below} $x=(x_1,x_2) \in V_k$ iff $x_1 = y_1$ in \attr and $x_2 > y_2$. In this case, one also says $x$ is \emph{directly above} $y$.
\end{defn}

\begin{lemma}\label{thm:directly-above}
  For any $\gx \in \V_k$, there is a vertex $(\gx,k+1)$ in $V_{k+1}$ directly above $(\gx,k)$.
  \begin{proof}
    If $\gx \in \V_k$, then $\gx = \simt_w(q_j)$ for some $w \in \Words_k$ and $j \in \{0,1,\dots,J\}$, and 
    \linenopax
    \begin{align*}
      (\gx,k)
      = (\simt_w(q_j),k) 
      \nbr (\simt_{wj} q_j,k+1) 
      = (\simt_{w} q_j,k+1)
      = (\gx,k+1)
    \end{align*}
    since $q_j$ is the fixed point of $\simt_j$. 
  \end{proof}
\end{lemma}

The previous lemma states that every vertex has a neighbour directly above itself. It is easy to see that a vertex $x = (\gx,k) \in V_\ell$ also has a neighbour directly below itself except in the case when $\gx \in \V_k \less \V_{k-1}$ and $k=\ell$.


%% file: random-walk-measures.tex

\section{The simple random walk on a network}
\label{sec:random-walks}

Here and elsewhere, the notation $V$ and $V_k$ continue to refer to the vertex set of \NF as in Definition~\ref{def:Vk-and-Gk}. 

\begin{defn}\label{def:Gamma(a,b)}
  An \emph{(infinite) path} in a network is an infinite sequence of vertices $\gw = (x_n)_{n \in \bN}$, where  $x_n \nbr x_{n+1}$ for all $n$. 
  Let \Paths be the space of all infinite paths \cpath in \NF. 
  Similarly, a \emph{finite path} $\gg = (y_n)_{n=0}^N$ is an finite sequence of vertices with $y_n \nbr y_{n+1}$ for all $n$. The \emph{length} of the finite path $\gg = (y_n)_{n=0}^N$ is denoted by $|\gg| = N$.
  
  A finite path $\gg  = (y_n)_{n=0}^N$ is also used to denote the cylinder set in $\Paths$ consisting of all infinite paths starting with the specified finite path:
  \linenopax
  \begin{align}\label{eqn:cylinder-set}
    \Paths(\gg) = \{\cpath=(x_n)_{n \in \bN} \in \Paths \suth x_n = y_n, n=0,1,\dots,|\gg|\}.
  \end{align}
  If $\gg = (x)$ is a path of length $0$, then $\Paths(\gg) = \Paths(x)$ is the collection of paths starting at $x$.
\end{defn}

\begin{defn}\label{def:prob}
  For any probability measure \gn on \NF, the space \Paths carries a natural probability measure $\prob_\gn$ for which 
  \linenopax
  \begin{align}\label{eqn:prob}
    \prob_\gn\left(\Paths(\gg)\right) = \gn(y_0) \prod_{n=1}^{N} p(y_{n-1},y_n),
  \end{align} 
  where $p(x,y) = c_{xy}/c(x)$ is the transition probability of the random walk; see \cite{Dynkin}. For $\gn=\gd_x$, we write $\prob_x := \prob_{\gd_x}$ for this probability measure on $\Paths(x)$. 
  We also use $\Ex_\gn$ to denote expectation with respect to $\prob_\gn$ and $\Ex_x[f] = \int f\,d\prob_x$.
\end{defn}

\begin{defn}\label{def:Xn}
  The \emph{random walk} started at a point $x_0$ on a network is the reversible Markov chain denoted by $(X_n)_{n=0}^\iy$, where $X_0 = x_0$ and $X_n$ is a $V$-valued random variable on \Paths for each $n=1,2,\dots$. Thus $X_n(\cpath)$ is the \nth coordinate of \cpath and denotes the location of the random walker at time $n$. As mentioned in Definition~\ref{def:prob}, the likelihood that a random walker at $x \in V$ steps to $y \nbr x$ is given by $p(x,y) = c_{xy}/c(x)$. More generally, the probability that a random walker at $x \in V$ will be at $z$ after taking $N$ steps is given by
  \linenopax
  \begin{align*}
    p_N(x,z) = \sum_{\gg \in \Paths_N(x;z)} \prod_{n=1}^{N} p(y_{n-1},y_n),
  \end{align*}
  where
  \linenopax
  \begin{align*}
    \Paths_N(x;z) 
    = \{\gg = (y_n)_{n=0}^N \suth |\gg|=N, y_0=x, y_N = z, \text{ and } y_n \nbr y_{n+1} \text{ for } n=1,\dots,N\}. 
  \end{align*}
  If $P$ is the infinite stochastic matrix with entries $[P]_{x,y} = p(x,y)$, then 
  $[P^N]_{x,y} = p_N(x,y)$. The \emph{Green function} is 
  \linenopax
  \begin{align}\label{eqn:Green-function}
    G(x,y) := \sum_{N=1}^\iy p_N(x,y).
  \end{align}
\end{defn}

$G(x,y)$ should not be confused with the Green function on \attr as discussed in \cite[\S3.5]{Kig01}, \cite[\S2.6]{Str06}, or \cite{resolvent}.

\begin{defn}\label{def:transient}
  A network is \emph{transient} if the simple random walk leaves any finite subset almost surely. In other words, if $A$ is a finite subnetwork, then there exists $N$ such that $\prob[X_n \in A \suth n \geq N] = 0$. The Green function converges for some choice of $x,y \in \NF$ (or equivalently, for all $x,y \in \NF$) if and only if the network is transient.
\end{defn}

The next result is superceded by Theorem~\ref{thm:Convergence-to-the-boundary}, but is heuristically helpful.

\begin{figure}
  \scalebox{0.90}{\includegraphics{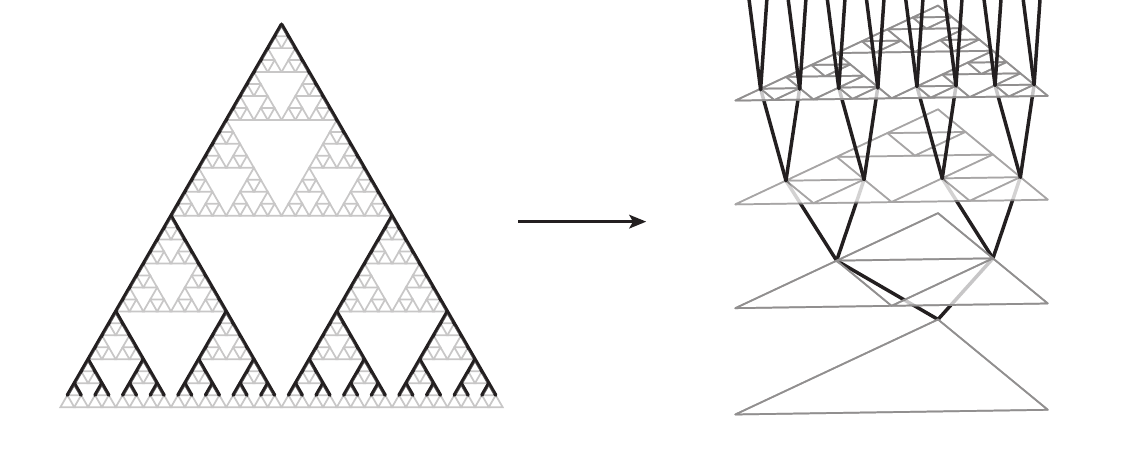}}
  \centering
  \caption{\captionsize An embedding of the binary tree into the Sierpinski network $\network_{\negsp[4]\sS\sG}$.}
  \label{fig:sierpinski-subtree}
\end{figure}

\begin{scholium}\label{thm:network-is-transient}
  The network \NF is transient.
  \begin{proof}[Sketch of proof]
    It is well-known that it suffices to find a transient subnetwork of \NF, see \cite[Cor.~2.15]{Woess}, for example. Figure~\ref{fig:sierpinski-subtree} depicts an embedding of the binary tree into $\network_{\negsp[4]\sS\sG}$, and it is well-known that the binary tree is transient.\footnote{An easy way to prove this is to note that the function $w(x) = 2^{-dist(x,o)}$ (where $dist(x,y)$ is the number of edges between $x$ and $y$) defines a flow of finite energy from $o$ to infinity, and apply \cite[Thm.~4.51]{WoessII} or \cite[Thm.~2.10]{Lyons}.} 
    Since the iterated function system always contains at least two maps $\simt_0$ and $\simt_1$, it is easy to see how one can embed the binary tree in a general network \NF in a similar fashion.
  \end{proof}
\end{scholium}

\begin{defn}\label{def:path-to-infinity}
  One calls $\cpath \in \Paths(o)$ a \emph{path to infinity} iff $\cpath = (x_n)_{n \in \bN}$ eventually leaves every finite subset of $V$. That is, for a finite set $A \ci V$, there must be some $N$ such that $x_n \notin A$ for all $n \geq N$. In this case, we write $\gw \to \iy$ or $x_n \to \iy$. 
  For a path $\cpath \in \Paths(o)$, the limiting value of $f$ along \cpath is
  \linenopax
  \begin{align}\label{eqn:path-limit-of-f}
    \lim\nolimits_\cpath f := \lim_{n \to \iy} f(x_n),
    \qq \text{for }\; \cpath  = (x_n)_{n \in \bN}.
  \end{align}
\end{defn}

\begin{remark}[``At infinity'']
  \label{rem:at-infinity}
  The imprecise terminology \emph{at infinity} is used in two different senses in this paper: one refers to vertices of \NF with second coordinate tending to \iy, and the other refers to the location $X_n(\gw)$ of the random walker, as $n \to \iy$. 
  However, 
   it is shown in Theorem~\ref{thm:Convergence-to-the-boundary} that \prob is supported on the paths to infinity. That is, the paths which do not eventually leave every finite subset form a set of \prob-measure $0$, and consequently these two notions of ``at infinity'' agree up to a set of paths of \prob-measure $0$. 
\end{remark}

%% file: separating-points.tex

\section{Separating points at infinity}
\label{sec:separating-points}

This section shows how to construct a function $u$ on \NF which has finite energy 
and which has certain prescribed limiting behavior. 
The following definition allows one to think of the limit properties of a function on \NF in terms of ``boundary values''.

\begin{defn}\label{def:limits-in-SG}
  Let \gx be a point of \attr. If $\{(x_n,z_n)\}_{n=1}^\iy$ is a sequence in $V \ci \attr \times \bN$, then one says this sequence \emph{tends to \gx} as $n \to \iy$ iff  $x_n \to \gx$ with respect to the resistance metric on \attr and $z_n \to \iy$.%
  \footnote{For now, it suffices to consider convergence the usual Euclidean metric on \bRd, as regularity (see Axiom~\ref{axm:pcf}) ensures these two metrics induce the same topology; see \cite[Thm.~3.3.4]{Kig01} for this fact and \cite{ERM,Kig03,Lyons} for more on resistance metric.} 
  If $\cpath=\{(x_n,z_n)\}_{n=1}^\iy$ is a path in \NF, write $\cpath\to\gx$ if the vertices in \cpath form a sequence tending to \gx. 
  We also extend Definition~\ref{def:directly-below} so that $x = (x_0,z_0) \in\NF$ is \emph{directly below} $\gx \in SG$ iff $\{(x_0,k)\}_{k=z_0}^\iy$ tends to \gx.
\end{defn}

\begin{remark}[Boundary heuristics]
  \label{rem:Boundary-heuristics}
  The notion of a sequence in \NF converging to a point of \attr in Definition~\ref{def:limits-in-SG} is a temporary heuristic. More rigourously, we compactify the network \NF as in Theorem~\ref{thm:compactification2} and obtain the boundary $\BF = \NB \less \NF$ as in Definition~\ref{def:compactification}. After the boundary is shown to be (homeomorphic to) \attr in Theorem~\ref{thm:B=F}, it will be clear that sequences which converge in the sense of Definition~\ref{def:limits-in-SG} correspond to sequences in \NB which converge to points of $\BF \cong \attr$ in the usual sense.
  
  Similarly, the notion of functions on \NF on separating points of \attr is presented just below in Definition~\ref{def:separate-points}, and this should also be considered a heuristic. Later, these functions will be continuously extended to functions on \NB, where they will separate points of $\BF \cong \attr$ in the usual sense.
\end{remark}

The goal of this section is to construct functions of finite energy on \NF that separate points of \attr in the following sense. 

\begin{defn}\label{def:separate-points}
  For distinct points $\gx, \gz \in \attr$, it follows from \eqref{eqn:quotient} that there are $w^{(\gx)}, w^{(\gz)} \in \Words$ such that $\lim_{n \to \iy} \simt_{w^{(\gx)}|n}(q) = \gx$ and $\lim_{n \to \iy} \simt_{w^{(\gz)}|n}(q) = \gz$, for any $q \in \attr$. 
  A function $u$ on \NF is said to \emph{separate the points} \gx and \gz iff for any $q \in \attr$, one has
  \linenopax
  \begin{align*}
    \lim_{n \to \iy} u\left(\simt_{w^{(\gx)}|n}(q),n\right) = \ga
    \q\text{and}\q
    \lim_{n \to \iy} u\left(\simt_{w^{(\gx)}|n}(q),n\right) = \gb, 
  \end{align*}
  with $\ga \neq \gb$. See Remark~\ref{rem:Boundary-heuristics}.
\end{defn}

In fact, the construction of ``harmonically generated function'' (see Definition~\ref{def:harmonically-generated-function}) presented below allows one to choose arbitrarily small disjoint neighbourhoods $U$ of \gx and $W$ of \gz ($U,W \ci \attr$), and then construct a function of finite energy $u$ on \NF for which $\lim_\cpath u > 0$ if \cpath tends to point of $U$, and $\lim_\cpath u < 0$ if \cpath tends to a point of $W$, and $\lim_\cpath u = 0$ if \cpath tends to some point of $\attr \less (U \cup W)$, whenever \cpath is a path in \NF that tends to some point of \attr in the sense of Definition~\ref{def:limits-in-SG}.

\subsection{Energy on the Sierpinski network}

By construction, 
\linenopax
\begin{align*}
  E \cup F = \bigcup_{k=0}^\iy E_k \;\cup\; \bigcup_{k=0}^\iy F_k
\end{align*}
is a partition of the edges of the Sierpinski network. The energy of a function $f:G \to \bR$ is a sum over the edges of the network and thus decomposes in terms of this partition as
\linenopax
\begin{align}\label{eqn:energy-decomp-by-levels}
  \energy(f) = \energy_E(f) + \energy_F(f)
  &= \sum_{k=0}^\iy \energy_{E_k}(f) + \sum_{k=1}^\iy \energy_{F_k}(f),
\end{align}
where $\energy_H(f) = \sum_{(xy) \in H} c_{xy} (f(x)-f(y))^2$ denotes the restriction of the energy to a subset $H$ of the edges of the graph. The sum $\energy_H$ is understood to contain each edge only once, so that exactly one term appears in the sum for each (undirected) edge $\{(xy), (yx)\}$ in $H$.

Considering $\simt_j$ as a map $\simt_j:\attr \to \attr$, define a mapping $\gY_j: \attr \times \bZ_+ \to \attr \times \bZ_+$ by 
\linenopax
  \begin{align}\label{eqn:psi_j}
  \gY_j(x,z) = (\simt_j(x),z+1),
  \qq x \in \attr, z \in \bZ_+, j=0,1,2, 
\end{align}
and let $\gY_w$ denote a composition of the mappings $\gY_j$ as in \eqref{eqn:word-composition}.

\begin{defn}[Harmonically generated function]
  \label{def:harmonically-generated-function}
  A \emph{harmonically generated function} $u$ on \NF is completely specified by its values on $V_0$. Let  
	\linenopax
  \begin{align}\label{eqn:u-initial-values}
	  u_0 = \left[\begin{array}{cccc}
	    u(q_0) & u(q_1) & \dots & u(q_J)
	  \end{array}
	  \right]^T
	\end{align} 
	be given. Then assign the values of $u$ on the vertices of $\gG_w = \gY_w \gG_0$ via
\linenopax
\begin{align}\label{eqn:harmonic-extension-algorithm}
  u|_{\gY_w V_0} = A_w u|_{V_0},
  \qq \text{where } \; A_w = A_{w_k} \dots A_{w_2} A_{w_1},
\end{align}
  and the matrices $A_0, A_1, \dots, A_J$ are the \emph{harmonic extension matrices}. The complete definition of these matrices is lengthy but can be found in \cite[(3.2.2)]{Kig01} or \cite[(4.2.17)--(4.2.18)]{Str06}; see Remark~\ref{rem:harmonic-extension-algorithm}.
\end{defn}

\begin{remark}\label{rem:hg-not-harmonic}
  The algorithm \eqref{eqn:harmonic-extension-algorithm} will \textbf{not} typically produce a harmonic function on the Sierpinski network, but it will produce a function with nontrivial harmonic component, as in Definition~\ref{def:harmonic-component}. The nomenclature is intended to evoke the fact that these functions are generated by the harmonic extension procedure.
\end{remark}

\begin{remark}\label{rem:harmonic-extension-algorithm}
  The harmonic extension algorithm appears as part of the proof of \cite[Thm.~3.2.4]{Kig01} and is laid out more explicitly in \cite[\S3.1]{Str06}; see also \cite[\S4.3]{Str06}. The harmonic extension matrices $A_j$ in \eqref{eqn:harmonic-extension-algorithm} are stochastic, and so have 1 as their largest eigenvalue. The second-largest eigenvalue of $A_j$ coincides with the renormalization factor $r_j$ appearing in Definition~\ref{def:Vk-and-Gk} and Definition~\ref{def:Self-similar-energy-form}, provided $\attr \less V_0$ is connected; see \cite{Strichartz99} for this fact, and \cite[A.1.2]{Kig01} for a more general result.
\end{remark}


With \gY replaced by \gF in \eqref{eqn:harmonic-extension-algorithm}, this procedure \eqref{eqn:harmonic-extension-algorithm} is called the \emph{harmonic extension algorithm} in the theory of analysis on fractals \cite{Kig01,Str06}. It provides for an explicit construction of a harmonic function on a pcf self-similar fractal with prescribed values on $V_0 = \{q_0,q_1,\dots,q_J\}$.%
  \footnote{These are called ``boundary values'' in \cite{Kig01,Str06}, but in a different sense than is considered in this paper.} 
  
For the fractal \SG, the harmonic extension matrices $A_j$ are given by
\linenopax
\begin{align*}
  A_0 = \left[
  \begin{array}{rrr}
    1 & 0 & 0 \vstr\\ 
    \frac25 & \frac25 & \frac15 \vstr\\ 
    \frac25 & \frac15 & \frac25 \vstr\\
  \end{array}\right],
  \qq
  A_1 = \left[
  \begin{array}{rrr}
    \frac25 & \frac25 & \frac15 \vstr\\ 
    0 & 1 & 0 \vstr \\ 
    \frac15 & \frac25 & \frac25 \vstr\\
  \end{array}\right],
  \qq
  A_2 = \left[
  \begin{array}{rrr}
    \frac25 & \frac15 & \frac25 \vstr\\ 
    \frac15 & \frac25 & \frac25 \vstr\\
    0 & 0 & 1 \vstr\\ 
  \end{array}\right],
\end{align*}
  see \cite[(1.3.28)]{Str06}, or \cite[Ex.~3.2.6]{Kig01}. See Figure~\ref{fig:sierpinski-bd-fn} for an illustration of the harmonically generated function on the Sierpinski network $\network_{\sS\sG}$ with initial values $u(q_o)=1, u(q_1)=0, u(q_2)=0$.  

\begin{lemma}\label{thm:harmext-is-finite-energy}
  For any choice of initial values $u_0 = [a_0,a_1,\dots,a_J]^T$ on $V_0$ as in \eqref{eqn:u-initial-values}, the harmonically generated function $u$ defined on \NF via \eqref{eqn:harmonic-extension-algorithm} has finite energy.
  \begin{proof}
    The harmonic extension algorithm has two useful features which facilitate the use of \eqref{eqn:energy-decomp-by-levels} in computing the energy of $u$. Following the notation of Definition~\ref{def:Vk-and-Gk}, let $\rmax := \max\{r_0,r_1,\dots,r_J\}$ be the largest of the renormalization factors, and observe that 
    \linenopax
    \begin{align}\label{eqn:k-level-energy}
      \energy_{E_k}(u) \leq \rmax \energy_{E_{k-1}}(u),  
    \end{align}
    Note that $0 < \rmax < 1$ by Axiom~\ref{axm:pcf}, which implies that $\sum_{k=0}^\iy \energy_{E_k}(u)$ is a convergent geometric series. 
  Second, consider $\energy_{F_k}(u)$. If we pick a vertex $x \in V_k$, then it follows from the harmonic extension procedure that $u(y) = u(x)$ for any vertex $y$ which is directly above $x$, and so there will be no contribution to the sum from any such edges. Neglecting these elements of $F_k$, Definition~\ref{def:vertical-edge} and the harmonic extension procedure imply that each nonzero summand in $\energy_{F_k}(u)$ appears also in $\energy_{E_k}(u)$, and hence $\energy_{F_k}(u) \leq \energy_{E_k}(u)$.%
  \footnote{For an illustration of this fact, compare the energy contributions from $F_k$ in Figure~\ref{fig:sierpinski-network-v-edges} to those of $E_k$ in Figure~\ref{fig:sierpinski-network-h-edges}.}
  Combining these facts with the decomposition \eqref{eqn:energy-decomp-by-levels} gives
  \linenopax
  \begin{align*}
      \energy(u) 
      &= \sum_{k=0}^\iy \energy_{E_k}(u) + \sum_{k=1}^\iy \energy_{F_k}(u) 
      \leq 2 \sum_{k=0}^\iy \rmax ^k \energy_{E_0}(u) 
      = \frac{\energy_{E_0}(u)}{1-\rmax},
    \end{align*}
    and shows that $\energy(u) < \iy$.
  \end{proof}
\end{lemma}

\subsection{Harmonically generated functions have a nontrivial harmonic component}

We will also need two results from the literature; one classical and one more recent. 

\begin{defn}\label{def:Dirichlet-norm}
  Define the \emph{Dirichlet norm} of $f:V \to \bR$ by 
  \linenopax
  \begin{align}\label{eqn:Dirichlet-norm}
    \energy_o(f) = \energy(f) + \cond(o) f(o)^2.
  \end{align}
\end{defn}

\begin{prop}[{\cite[Lem.~1.3]{Yamasaki79}}]
  \label{thm:Royden-decomp}
  The space of Dirichlet functions $\Gdd = \{f \suth \energy_o(f) < \iy\}$ can be written as $\Gdd = \Gddo + \GHD$, where \Gddo consists of those elements of \Gdd which are $\energy_o$-limits of functions of finite support, and \GHD consists of the elements of \Gdd which are harmonic everywhere.
\end{prop}

The representation of \Gdd given by Proposition~\ref{thm:Royden-decomp} is often called the (discrete) Royden decomposition, in reference to H. Royden's analogous result for Riemannian manifolds. 

\begin{defn}\label{def:harmonic-component}
  For $f \in \Gdd$, there corresponds a unique decomposition $f = f_{\Gddo} + f_{\GHD}$. We say that $f_{\GHD}$ is the \emph{harmonic component} of $f$.
\end{defn}

\begin{theorem}[{\cite[Cor.~1.2]{AnconaLyonsPeres}}]
  \label{thm:nonvanishing-Dirichlet-limits}
  If $X_n$ is the simple random walk (as in Definition~\ref{def:Xn}) on a transient network and $f = f_{\Gddo} + f_{\GHD}$ is the Royden decomposition of $f \in \Gdd$, then $\lim f(X_n) = \lim f_{\GHD}(X_n)$ \prob-a.s., and the limit is a random variable in $L^2(\Paths,\prob)$.
\end{theorem}

In other words, Theorem~\ref{thm:nonvanishing-Dirichlet-limits} states that with respect to the usual path space measure \prob (see Definition~\ref{def:prob}), $f_{\Gddo}$ tends to $0$ along almost every infinite path (and hence on almost every path to infinity). 

\begin{lemma}\label{thm:u-nonvanishing-at-infinity}
  Let $u_0 = [1,0,\dots,0]^T$ and construct $u$ via the harmonic extension algorithm as in Lemma~\ref{thm:harmext-is-finite-energy}. Then $u$ extends to a continuous function on the compactification $\NF \cup \attr$.
  \begin{proof}
    If \gx is a fixed rational point of \attr (see Definition~\ref{def:junction-point}) and $x_n = (\gx,n)$, $n \in \{n_0, n_0+1, n_0+2, \dots\}$, then $\{x_n\}_{n=n_0}^\iy$ is a sequence in \NF tending to \gx. 
  \end{proof}
\end{lemma}

\begin{lemma}\label{thm:u-has-harmonic}
  Let $u_0 = [1,0,\dots,0]^T$ and construct $u$ via the harmonic extension algorithm as in Lemma~\ref{thm:harmext-is-finite-energy}. Then $u$ has a nontrivial harmonic component, that is, $u_{\GHD} \neq 0$. 
  \begin{proof}
    Recall that $V_0 = \{q_0,q_1,\dots,q_J\}$.  
    From Definition~\ref{def:harmonically-generated-function}, it follows that $u(x) > 0$ unless $x$ is directly below an element of $\{q_1,\dots,q_J\}$, in which case $u(x)=0$. 
    The harmonic extension algorithm and Definition~\ref{def:harmonically-generated-function} imply that 
    \linenopax
    \begin{align*}
      \begin{cases}
        \lim_{n \to \iy} u(x_n) = 0, &\gx \in \{q_1,\dots,q_J\}, \\
        \lim_{n \to \iy} u(x_n) > 0, & \text{otherwise}.
      \end{cases}
    \end{align*}
    It is clear by symmetry that the set of paths in $\Paths(o)$ tending to $\{q_1,\dots,q_J\}$ is a set of \prob-measure $0$, so Theorem~\ref{thm:nonvanishing-Dirichlet-limits} implies that $u_{\GHD} \neq 0$.
  \end{proof}
\end{lemma}

\subsection{``Localizing'' the construction of $u$} 

In this section, we modify the function $u$ constructed in Lemma~\ref{thm:u-nonvanishing-at-infinity} so that the corresponding harmonic component $u_{\GHD}$ vanishes at infinity, except on paths tending to a chosen neighbourhood of a selected point in \attr.

\begin{defn}\label{def:boundary-point-of-an-m-cell}
  A \emph{boundary point of an $m$-cell} $K \ci \attr$ is a point $\gx \in K$ for which $\gx \in U \ci \attr$ implies $U \not \ci K$, whenever $U$ is open in the metric topology of effective resistance (or equivalently, in the subspace topology inherited from \bRd). The finite set of boundary points of a cell $K$ will be denoted by $\del K$.
\end{defn}

\begin{defn}\label{def:x-cell}
  Denote by $K(\gx,m)$ the union of the $m$-cells of \attr containing \gx.
\end{defn}

If \gx is a junction point, the set $K(\gx,m)$ consists of at least two (but at most finitely many) $m$-cells. If \gx is an irrational point or a rational point which is not a junction point, then $K(\gx,m)$ consists of just one $m$-cell; see Definition~\ref{def:junction-point}.

\begin{defn}[(Localized) harmonically generated functions]
  \label{def:localized-harmonically-gen'd}
  Recall from Definition~\ref{def:junction-point} that every rational point $\gx \in \attr$ has a generation of birth $m_\gx$; see \eqref{eqn:generation}. For $m \geq m_\gx$, the vertex $(\gx,m)$ is the unique point of $\gG_m$ which is directly below \gx, and \gx is a boundary point of an $m_\gx$-cell in \attr.  

For a fixed $x = (\gx,m) \in \NF$, define the \emph{(localized) harmonically generated function} $u_x$ by repeating the construction of $u$ via the harmonic extension procedure of Lemma~\ref{thm:harmext-is-finite-energy} on the set $\NF \cap \left(K(\gx,m) \times [m,\iy)\right)$ and extending $u_x$ by $0$ elsewhere.
 
More precisely, let $Q$ be the largest subset of $V_0$ such that for each $q \in Q$ there is a $w^{(q)} \in \Words_m$ with $\gY_{w^{(q)}}(q) = \gx$, and define 
\linenopax
\begin{align}\label{eqn:u^K}
  u_x(y) 
  := \begin{cases}
    u(\gY_{w^{(q)}}^{-1}(y)), &x \in \gY_{w^{(q)}}(\NF) \text{ for some } q \in Q,\\
    0 &\text{else}.
  \end{cases}
\end{align}
  Note that $\gY_{w^{(q)}}(\NF)$ and $\gY_{w^{(q')}}(\NF)$ are disjoint except for the vertices directly below \gx, for $q,q' \in Q$ with $q \neq q'$. For either $q$ or $q'$, formula \eqref{eqn:u^K} indicates that $u_x$ takes the value $1$ on the set of vertices directly below \gx, and hence \eqref{eqn:u^K} is well-defined and unambiguous.
\end{defn}

\begin{lemma}\label{thm:E(u_x)-is-finite}
  For $x=(\gx,m) \in \NF$, the harmonically generated function $u_x$ of Definition~\ref{def:localized-harmonically-gen'd} has finite energy.
  \begin{proof}
    Just as in the proof of Lemma~\ref{thm:u-nonvanishing-at-infinity}, one can see that $u_x$ vanishes on the set of vertices directly below the boundary points of $K(\gx,m)$, viz., the neighbours of $\gx^{(m)}$ in $\gG_m$, and the vertices directly above them. Since $u_x$ is extended by $0$ on the complement of $\bigcup_{q \in Q} \gY_{w^{(q)}}(\NF)$, there is no contribution to the energy sum \eqref{eqn:energy-decomp-by-levels} from edges exterior to $\bigcup_{q \in Q} \gY_{w^{(q)}}(\NF)$, except for the edges connecting $\gx^{(k)}$ to its neighbours in $\gG_{k-1}$. Denote this contribution to the energy of $u_x$ by $E$. Then $\energy(u_x) \leq C\energy(u) + E$, where $C=1$ if \gx is a rational point which is not a junction point, and $C=N$ if \gx is a junction point at which $N$ $m$-cells intersect. The rest of the argument for $\energy(u_x) < \iy$ and $u_x \neq 0$ is directly analogous to the proof of Lemma~\ref{thm:u-nonvanishing-at-infinity}.
  \end{proof}
\end{lemma}

\begin{theorem}[Separation]\label{thm:u_x-in-domE}
  Any two points $\ga,\gb \in \attr$ can be separated by a harmonically generated function $u_x$. 
  \begin{proof}
    If $\ga,\gb \in \attr$ are distinct, one can always find $m_\ga, m_\gb \in \bZ_+$ large enough to ensure $K(\ga,m_\ga) \cap K(\gb,m_\gb) = \es$ because the sets $\inter K(\ga,m)$ form a neighbourhood basis at \ga in the metric topology of \attr (which is obviously Hausdorff), and similarly for \gb. 
    Therefore, one can choose $x = (\gx,m) \in \NF$, such that $\ga \in \inter(K(\gx,m))$ and $\gb \notin K(\gx,m)$. Now $\lim_\gw u_x > 0$ for any \gw tending to \ga, and $\lim_\gg u_x = 0$ for any \gg tending to \gb.
  \end{proof}
\end{theorem}

\begin{cor}\label{thm:nontrivial-harmonics}
  Each harmonically generated function $u_x$ has a nontrivial harmonic component; in the notation of Definition~\ref{def:harmonic-component}, this means $(u_x)_{\GHD}$ is a bounded and (globally) harmonic function on \NF, and $\energy((u_x)_{\GHD}) < \iy$.
  \begin{proof}[Sketch of proof]
    Boundedness follows from the fact that $\sup_{y \in \NF} |u_x(y)| = 1$ and $u_x \in \dom \energy$ imply that $\sup_{y \in \NF} |(u_x)_{\GHD}(y)| < \iy$; see \cite{Soardi94, bdG, DGG}, for example.
    For the other claim, let $x = (\gx,0) \in V_0$, and suppose $\gw \to \iy$. Since $\lim_\gw u_x > 0$ unless \gw tends to a point $\gz \in \V_0 \less\{\gx\}$, it is clear from Theorem~\ref{thm:u_x-in-domE} (and by symmetry) that $u_x$ is nonvanishing at \iy on a set of positive $\Prob_x$-measure.%
    \footnote{In fact, it follows from Theorem~\ref{thm:Convergence-to-the-boundary} that $u_x$ is supported at \iy $\Prob_x$-a.s. The details of this proof can readily be filled in after reading \S\ref{sec:convergence}.} 
    This argument extends readily to other $x \in \NF$. 
    The result then follows from Theorem~\ref{thm:nonvanishing-Dirichlet-limits}.
  \end{proof}
\end{cor}

\begin{remark}\label{rem:nontrivial-harmonics}
  Consider the collection
  \linenopax
  \begin{align}\label{eqn:xF}
    \xF = \{u_x \suth x \in \NF\}.
  \end{align}
  Since each $u_x$ is bounded and \bR-valued by construction, it follows that the boundary constructed by applying Theorem~\ref{thm:compactification} to \eqref{eqn:xF}  contains (a homeomorphic image of) \attr. In Theorem~\ref{thm:B=F}, we will use Theorem~\ref{thm:compactification2} to prove the opposite containment. This requires development of the path-space formalism in \S\ref{sec:convergence}.
\end{remark}

%% file: convergence.tex

\section{Convergence to the boundary}
\label{sec:convergence}

The goal of this section is to establish the convergence of the random walk on \NF to the boundary, in precise sense described in Theorem~\ref{thm:Convergence-to-the-boundary}. This section follows the remarkably lucid approach of \cite[Ch.~7]{WoessII} and \cite{Dynkin}.

Proposition~\ref{thm:sequence-crossings} will be used in conjunction with the beautiful crossings estimate of \cite{AnconaLyonsPeres} in Theorem~\ref{thm:crossings-estimate} to obtain convergence results in the proof of Theorem~\ref{thm:Convergence-to-the-boundary} (actually, in Lemma~\ref{thm:prob_x(Paths_iy) = 1}). Here and elsewhere, the notation $C\left((r_n)_{n=1}^\iy \vr [a,b]\right)$ indicates the number of crossings of the interval $[a,b]$ by the sequence $(r_n)_{n=1}^\iy \ci \bR$. 

\begin{prop}\label{thm:sequence-crossings}
  For a sequence $(r_n)_{n=1}^\iy$ in \bR, $\lim r_n$ exists in $[-\iy,\iy]$ if and only if    
  \linenopax
  \begin{align}\label{eqn:sequence-crossings}
    C\left((r_n)_{n=1}^\iy \vr [a,b]\right) < \iy, 
    \qq \text{ for every interval } [a,b].
  \end{align}
  Without loss of generality, one may restrict to intervals with $a,b \in \bQ$ and $a<b$.
\end{prop}

\begin{theorem}[{\cite[Thm.~1.3]{AnconaLyonsPeres}}]
  \label{thm:crossings-estimate}
  Let $X_n$ be the random walk on \NF, and let $f:\NF \to \bR$. Then for any $a<b$,
  \linenopax
  \begin{align}\label{eqn:crossings-estimate}
    \Ex_x\left[C \left(f(X_n) \vr [a,b] \right) \right] 
    \leq C_x \frac{\energy(f)}{(b-a)^2},
  \end{align}
  where the constant $C_x := 2G(x,x)/c(x)$ depends only on $x$ and $G$ is the Green function. 
\end{theorem}

\begin{defn}\label{def:boundary}
  Define an equivalence relation on \Paths by 
  \linenopax
  \begin{align}\label{eqn:boundary-relation}
    \gw \sim \gg \q\iff\q \lim_\gw f = \lim_\gg f, 
    \qq\text{for all } f \in \xF.
  \end{align}
  The \emph{boundary} of \NF is then
  \linenopax
  \begin{align}\label{eqn:boundary}
    \BF := \Paths / \sim,
  \end{align}
  and $\NB = \NF \cup \BF$ is a compactification of \NF by Theorem~\ref{thm:compactification2}. Switching to the alternative formulation of Theorem~\ref{thm:compactification}, for any $f \in \xF$ and $x_\iy \in \BF$ with representative $\gw = (x_n)$, one can define
  \linenopax
  \begin{align}\label{eqn:f(gx)}
    f(x_\iy) := \lim_\gw f.
  \end{align}
  We will be particular interested in the set $\Paths_\iy \ci \Paths$ defined by 
  \linenopax
  \begin{align}\label{eqn:Paths_infty}
    \Paths_\iy := \{(x_n)_{n=1}^\iy \in \Paths \suth x_n \to x_\iy \in \BF
    \text{ in the topology of \NB}\}.
  \end{align}
\end{defn}

\begin{lemma}\label{thm:Paths_infty-is-measurable}
  The set $\Paths_\iy$ is measurable with respect to \sA, the Borel \gs-algebra of \NB.
  \begin{proof}
    We will boil down $\Paths_\iy$ in terms of cylinder sets. 
    To begin, note that $\Paths_\iy = \bigcap_{x \in \NF} \Paths_x$, where
    \linenopax
    \begin{align}\label{eqn:Paths_x}
      \Paths_x := \{\gw \in \Paths \suth \lim\nolimits_\gw u_x
      \text{ exists in \bR}\}.
    \end{align}
    Since 
    $\left(u_x(x_n)\right)_{n=1}^\iy$ is a bounded sequence, one can write 
    \linenopax
    \begin{align*}
      \Paths_x = \bigcap_{(a,b) \in \bQ^2} A_x([a,b]),
    \end{align*}
    where
    \linenopax
    \begin{align*}
      A_x([a,b]) := \{(x_n) \suth C\left((u_x(x_n)) \vr [a,b]\right) < \iy\}. 
    \end{align*}
    Now one can check that $\Paths_\iy \less A_x([a,b]) \in \sA$ for each fixed $[a,b]$, since
    \linenopax
    \begin{align*}
      \{(x_n) \suth C\left((u_x(x_n)) \vr [a,b]\right) < \iy\}
      = \bigcap_{m \geq 1} \bigcup_{j,k \geq m} 
      \left(\{(x_n) \suth u_x(x_j) \geq b\} \cup \{(x_n) \suth u_x(x_k) \geq b\}\right).
    \end{align*}
    Finally, 
    \linenopax
    \begin{align*}
      \{(x_n) \suth u_x(y_j) \geq b\} 
      = \bigcup_{u_x(x_j) \geq b} \Paths(y_0,\dots,y_j),
    \end{align*}
    where $\Paths(y_0,\dots,y_l)$ is a cylinder set as in \eqref{eqn:cylinder-set}, and similarly for 
    $u_x(y_k) \leq a$. 
  \end{proof}
\end{lemma}

\begin{lemma}\label{thm:prob_x(Paths_iy) = 1}
  For every $x \in \NF$, one has $\prob_x(\Paths_\iy) = 1$. 
  \begin{proof}
    We must show that $\lim_{n \to \iy} f(X_n)$ exists $\prob_x$-almost surely, for every function in the family \xF of \eqref{eqn:xF}. From Theorem~\ref{thm:crossings-estimate}, we have
    \linenopax
    \begin{align*}
      \Ex_x\left[C \left(u_x(X_n) \vr [a,b] \right) \right] 
      \leq 2 \frac{G(x,x)}{c(x)} \frac{\energy(u_x)}{(b-a)^2}
       < \iy,
    \end{align*}
    whence for every fixed interval $[a,b]$ one has
    \linenopax
    \begin{align*}
      C\left(u_x(X_n))_{n=1}^\iy \vr [a,b]\right) < \iy, 
      \qq \prob_x\text{-a.s.}, 
    \end{align*}
    and therefore $\lim_{n \to \iy} u_x(X_n)$ exists $\Pr_x$-almost surely. 
  \end{proof}
\end{lemma}

\begin{lemma}\label{thm:prob_x(Paths_iy) = 1}
  For each $x \in \NF$, the function $X_\iy:\Paths \to \BF$ defined by     
  \linenopax
  \begin{align}\label{eqn:X_infty}
    X_\iy(\gw) = x_\iy, \qq\text{for all } \gw = (x_n) \in \Paths_\iy
  \end{align}
  is measurable with respect to \sA.
  \begin{proof}
    In view of \eqref{eqn:xF} and \eqref{eqn:boundary-relation}, the topology on \BF has a subbasis consisting of the sets
    \linenopax
    \begin{align}\label{eqn:subbasicU}
      U_{x,x_\iy,\ge} &:= \{y_\iy \in \BF \suth |u_x(x_\iy) - u_x(y_\iy)|<\ge\},
    \end{align}
    for any $x \in \NF$, $x_\iy \in \BF$, and $\ge>0$. (Note that $x_\iy$ is not related to $x$ in any way; they are free to vary separately.)
    Observe from \eqref{eqn:f(gx)} that $u_x$ 
    is defined at $x_\iy,y_\iy \in \BF$ by continuity. Since
    \linenopax
    \begin{align*}
      [X_\iy \in U_{x,x_\iy,\ge}]
      = \{(y_n) \in \Paths_\iy \suth |u_x(x_\iy) - \lim_{n \to \iy} u_x(y_n)|<\ge\},
    \end{align*}
    and similarly for $W_{x,x_\iy,\ge}$, one can see that $[X_\iy \in U_{x,x_\iy,\ge}]$, $[X_\iy \in W_{x,x_\iy,\ge}] \in \sA$.
  \end{proof}
\end{lemma}

\begin{theorem}[Convergence to the boundary]
  \label{thm:Convergence-to-the-boundary}
  There is a random variable $X_\iy$ such that for every $x \in \NF$,
  \linenopax
  \begin{align}\label{eqn:Convergence-to-the-boundary}
    \lim_{n \to \iy} X_n = X_\iy,
    \qq \prob_x-\text{almost surely},
  \end{align}
  in the topology of \NB.
  \begin{proof}
    This is a restatement of the combined results of Lemma~\ref{thm:Paths_infty-is-measurable}, Lemma~\ref{thm:prob_x(Paths_iy) = 1}, and Lemma~\ref{thm:Paths_infty-is-measurable}.
  \end{proof}
\end{theorem}

%% file: identification.tex

\section{Identification of the boundary}
\label{sec:identification}

In this section, we show that the boundary \BF of \NB (see Definition~\ref{def:boundary}) is homeomorphic to the fractal \attr. In reference to its extension to \BF by continuity, we think of the harmonically generated function $u_x$ (see Definition~\ref{def:localized-harmonically-gen'd}) as a ``bump function'' centered $\gx \in \attr$. More precisely, from the construction of $u_x$ and the extension by continuity in \eqref{eqn:f(gx)}, one has
    \linenopax
    \begin{align}\label{eqn:u=spline}
      \lim_\gw u_x = \lim_\gw u_{(\gx,m)} = \gy_\gx^{(m)}(\gz),
      \qq\text{ for any } \gw \to \gz \text{ and } x=(\gx,m),
    \end{align}
    where $\gy_\gx^{(m)}$ is the piecewise harmonic spline satisfying $\gy_\gx^{(m)}(y) = \gd_{xy}$ on $\V_m$; see \cite[Thm.~2.1.2 and \S2.2]{Str06}. 

\begin{defn}\label{def:u_K}
  For a cell $K \ci \attr$, define the ``$m$-cell of \NF at level $m$'' by
  \linenopax
  \begin{align}\label{eqn:K^(m)}
    K^{(m)} := \{x=(\gx,m) \suth \gx \in K\} = (K \times \{m\}) \cap \NF,
  \end{align}
  so that $K^{(m)}$ consists of the points of $\gG_m$ which lie directly below $K$, and define
  \linenopax
  \begin{align}\label{eqn:u_K^(m)}
    \uKm := \sum_{x \in K^{(m)}} u_x.
  \end{align}
  It is evident from Lemma~\ref{thm:u_K}(ii) that one can restrict the support of \uKm on \attr to an arbitrarily small neighbourhood of $K$ by choosing $m$ sufficiently large.
\end{defn}

\begin{lemma}\label{thm:u_K}
  For any cell $K \ci \attr$, the function \uKm is a bump function on $K$ with the following properties:
  \begin{enumerate}[(i)]
    \item For all $x \in K^{(l)}$ and $l \geq m$, one has $\uKm[l](x) = 1$. 
    \item For all points of \attr outside of $\bigcup_{x \in K^{(m)}} K(\gx,m)$, and all points of \NF directly below this set, the function $u_{K^{(m)}}$ vanishes.
  \end{enumerate}
  \begin{proof}
     The harmonic extension algorithm is a linear operation. Consequently, 
     \linenopax
     \begin{align}
       \gy_K := \sum_{\gx \in \del K} \gy_\gx^{(m)},
       \qq\text{where $K$ is an $m$-cell}
     \end{align}
     can be equivalently be constructed by summing the $\gy_\gx^{(m)}$, or by taking the piecewise harmonic spline on \attr which takes value 1 on $\del K$ and value 0 on the other rational points of generation $m$. From the latter formulation, it is clear that $\restr{\gy_K}{K} \equiv 1$. The same reasoning can be applied to verify claim (i). Claim (ii) is immediate from the construction of the $u_x$ in Definition~\ref{def:localized-harmonically-gen'd}. 
  \end{proof}
\end{lemma}

\begin{lemma}\label{thm:nice-reps}
  For $x_\iy \in \gW_\iy$, there is a representative $(x_n)_{n=0}^\iy$ of $x_\iy$ which can be written $(\gx_n,n)_{n=0}^\iy$. This representative can be chosen so that whenever $\gx_n$ lies in an $m$-cell $K$, then $\gx_{n+1} \in K$ also.
  \begin{proof}
    First, we show that for any representative $(x_n)_{n=0}^\iy = (\gx_n,k_n)_{n=0}^\iy$ and any fixed $m \in \bZ_+$, there is an $m$-neighbourhood $K(\gz,m)$ of $m$-cells containing the tail of $(\gx_n)_{n=0}^\iy$. Suppose this is not the case. Then there must be at least two $m$-neighbourhoods $K(\gz_1,m)$ and $K(\gz_2,m)$ such that $(\gx_n)$ is in each one for infinitely many $n$. In fact, the tail must be in each of $K(\gz_1,m) \less K(\gz_2,m)$ and $K(\gz_2,m) \less K(\gz_1,m)$ infinitely often, or else the tail is contained in one of the $m$-neighbourhoods. If there are two boundary points $\gh, \gw \in \del K(\gz_1,m)$ through which $(\gx_n)$ passes infinitely many times, then the limits of $u_{(\gh,m+1)}(x_n)$ and $u_{(\gw,m+1)}(x_n)$ cannot exist, which contradicts $x_\iy \in \gW_\iy$. So it would have to be the case that $(\gx_n)$ passes through a single boundary point $\gh \in \del K(\gz_1,m)$ infinitely many times. It must be the case for at least one of $i=1,2$ that $(\gx_n)$ leaves the $m$-cell of $K(\gx_i,m)$ containing \gh infinitely many times; otherwise $K(\gh,m)$ would contain the tail. Let \gw be the boundary point of $K(\gx_i,m)$ through which $(\gx_n)$ passes infinitely often. Then again the limits of $u_{(\gh,m+1)}(x_n)$ and $u_{(\gw,m+1)}(x_n)$ cannot exist, and now we are out of options.
  \end{proof}
\end{lemma}

\begin{theorem}\label{thm:B=F}
  The network boundary \BF is (homeomorphic to) the attractor \attr.
  \begin{proof}
    For $\gx \in \attr$, one has $\gx = \lim_{n \to \iy} \simt_{w|n}(q)$. We show that
    \linenopax
    \begin{align}\label{eqn:f-homeo}
      f:\attr \to \BF
      \qq\text{by}\qq 
      f:\gx \mapsto (\simt_{w|n}(q), n)_{n=0}^\iy 
    \end{align}
    defines a homeomorphism from $\attr$ to \BF. For $x_\iy \in \gW_\iy$, choose a representative $(x_n)_{n=0}^\iy$ according to Lemma~\ref{thm:nice-reps}, and denote $x_n = (\gx_n,n)$. Then $\gx_n$ is a rational point of \attr, so there is a $q \in \V_0$ and $w_1\dots w_n \in \Words_n$ such that $\gx_n = \simt_{w_1\dots w_n}(q)$. In fact, the nesting property in Lemma~\ref{thm:nice-reps} ensures there is an infinite word $w \in \Words$ for which $\gx_n = \simt_{w|n}(q)$, for every $n \in \bZ_+$.
    Now define 
    \linenopax
    \begin{align}\label{eqn:f-homeo}
      g:\BF \to \attr
      \qq\text{by}\qq 
      g:(\gx_n,n) \mapsto \lim_{n \to \iy} \simt_{w|n}(q). 
    \end{align}
    It is easy to verify that $g$ is the inverse of $f$.

    To see that $f$ is continuous, we show that the preimages of the subbasic open sets $U_{x,x_\iy,\ge}$ of \eqref{eqn:subbasicU} are open. 
    Recall from \eqref{eqn:u=spline} that $\restr{u_x}{\BF} = \gy_\gx^{(m)}$ for $x = (\gx,m)$.
    Let us use $\cj{\gx}$ to denote the representative of $f(\gx)$ chosen according to Lemma~\ref{thm:nice-reps}, so that $f(\gx) = [\cj{\gx}]$.
    \linenopax
    \begin{align}
      f^{-1}(U_{x,x_\iy,\ge})
      &= f^{-1}\left(\{ \cj{\gz} \in \BF \suth |u_x(\cj{\gx})-u_x(\cj{\gz})| < \ge\}\right) \notag \\
      &= \{\gz \in \attr \suth |\gy_x^{(m)}(\gx)-\gy_x^{(m)}(\gz)| < \ge\} \notag \\
      &= \left(\gy_x^{(m)}\right)^{-1} \left(B(\gy_x^{(m)}(\gx),\ge)\right).
      \label{eqn:preimage-of-subbasic}
    \end{align}
    It follows from \cite[(1.4.3)]{Str06} that $\gy_x^{(m)}$ is uniformly continuous on \attr, and hence \eqref{eqn:preimage-of-subbasic} is certainly open.
    We now have a continuous bijection from a compact space to a Hausdorff space, and hence a homeomorphism.
  \end{proof}
\end{theorem}